\documentclass[11pt]{article}
\usepackage[font=small,labelfont=bf]{caption}
\usepackage{latexsym,amsfonts,amsthm,amsmath,amssymb}
\usepackage{hyperref}
\usepackage{xcolor}
\usepackage{verbatim}
\usepackage{enumerate}
\usepackage{tabularx,graphicx}

\newtheorem{Def}{Definition}

\begin{document}

\begin{center}
{\LARGE Euler and the Legendre Polynomials} \\

\bigskip\bigskip

{Alexander Aycock, Johannes-Gutenberg University Mainz \\ Staudinger Weg 9, 55128 Mainz \\ \url{aaycock@students.uni-mainz.de}}

\bigskip

%{Maria Agnesi, University of Bologna \\ 86 Via del Pallone, Bologna, Italy 40100 \\ {\tt magnesi@bologna.edu}}

\end{center}

\medskip

%\title{Euler and the Legendre Polynomials}
%\author{Alexander Aycock, Johannes-Gutenberg University Mainz \\ Staudinger Weg 9, 55128 Mainz \\ \url{aaycock@students.uni-mainz.de} }
%\date{ }
%\maketitle

\begin{abstract} \noindent % Abstract not required for review articles.
In this note we will present how  Euler's investigations on various different subjects lead to certain properties of the Legendre polynomials. More precisely, we will show that the generating function and the difference equation for the Legendre polynomials was already written down by Euler in at least two different papers. Furthermore, we will demonstrate that some  familiar expressions for the Legendre polynomials  are corollaries of the before-mentioned works. Finally, we will show that Euler's ideas on continued fractions lead to an integral representation for the Legendre polynomials that seems to be less generally known. 
\end{abstract}

\section{Introduction}
\label{sec: Introduction}

The Legendre polynomials occur frequently at various places in physics. They were introduced by Legendre (1752-1833) in his paper  {\it ``Recherches sur l'attraction des sphéroïdes homogènes"} \cite{Le85} (``Researches on the attraction of homogeneous spheroids") in his study of the gravitational potential expressible as:

\begin{equation*}
    %\label{eq: Legendre and his Polynomials}
    \dfrac{1}{|\bf{x}-\bf{x}^{\prime}|}= \dfrac{1}{\sqrt{r^2 -2rr^{\prime} \cos (\beta)+(r^{\prime})^2}}= \sum_{n=0}^{\infty} \dfrac{{r^{\prime}}^n}{r^{n+1}}P_n(\cos (\beta)),
\end{equation*}
where $r$ and $r^{\prime}$ are the lengths of the vectors $\bf{x}$ and $\bf{x}^{\prime}$, respectively, $\beta$ is the angle between those two vectors, and $P_n(t)$ denotes a polynomial of degree $n$ in $t$.
They are also interesting from a mathematical viewpoint, because they exhibit a system of complete and orthogonal polynomials.

Nowadays, the $n$-th Legendre polynomial, $P_n(t)$, is often introduced as the coefficient a Taylor series. We have: 

\begin{equation}
    \label{eq: Definition Legendre polynomials}
    \dfrac{1}{\sqrt{1-2xt+x^2}}:= \sum_{n=0}^{\infty}P_n(t)x^n.
\end{equation}
Aside from Legendre's original definition and the more modern (\ref{eq: Definition Legendre polynomials}), there are many ways to define the Legendre polynomials and hence as many different explicit formulas for them. In this note we will first use their definition as solutions of the following difference equation\footnote{Later, in \ref{sec: Connection of this Result to similar Insights}, we will also use their definition as coefficients of a generating function.}:

\begin{equation}
\label{eq: Difference Equation for Legendre}
(n+1)P_{n+1}(t)= (2n+1)tP_n(t)-nP_{n-1}(t),\quad \text{for} \quad n \in \mathbb{N}_0.
\end{equation}
Requiring  that $P_0(t)=1$ and $P_1(t)=t$ this equation defines all Legendre polynomials uniquely. The first Legendre polynomials are seen in Table \ref{Table 1} below.

\begin{table}[h]
    \centering
    \begin{tabular}{c}
    \renewcommand{\arraystretch}{2,0}
\setlength{\arraycolsep}{1.5mm}
$\begin{array}{|l|l|} \hline
\quad  n \quad & \quad P_n(t) \\ \hline
\quad 0       & \quad 1 \\
\quad 1       & \quad t \\
\quad 2       & \quad \frac{1}{2}(3t^2-1) \\
\quad 3       & \quad \frac{1}{2}(5t^3-3t) \\
\quad 4       & \quad \frac{1}{8}(35t^4-30t^2+3) \\
\quad 5       & \quad \frac{1}{8}(63t^5-70t^3+15t) \\
\quad 6       & \quad \frac{1}{16}(231t^6-315t^4+105t^4-5) \\
\quad 7       & \quad \frac{1}{16}(429t^7-693t^5+315t^3-35t) \\ \hline
\end{array}$
    \end{tabular}
    \caption{The first Legendre polynomials for small values of $n$.}
    \label{Table 1}
\end{table}
We will show that (\ref{eq: Difference Equation for Legendre}) occurs as a special case of more general formulas in some  papers written by Euler (1707-1783). Those papers include {\it ``Speculationes super formula integrali $\int \frac{x^ndx}{\sqrt{aa-2bx+cxx}}$, ubi simul egregiae observationes circa fractiones continuas occurrunt"} \cite{E606} (E606: ``Speculations about the integral formula $\int \frac{x^ndx}{\sqrt{aa-2bx+cxx}}$, where at the same time extraordinary observations on continued fractions occur"), {\it ``Theorema maxime memorabile circa formulam integralem $\int \frac{\partial \phi \cos \lambda \phi}{(1+aa-2a\cos \phi)^{n+1}}$"} \cite{E672} (E672: ``A most memorable theorem on the integral formula $\int \frac{\partial \phi \cos \lambda \phi}{(1+aa-2a\cos \phi)^{n+1}}$"), {\it ``Disquitio coniecturalis super formula integrali $\int \frac{\partial \phi \cos i \phi}{(\alpha +\beta \cos \phi)^n}$"} \cite{E673} (E673, ``A conjectural investiagtion on the integral formula $\int \frac{\partial \phi \cos i \phi}{(\alpha +\beta \cos \phi)^n}$"), {\it ``Demonstratio theorematis insignis per coniecturam eruti circa intagrationem formulae $\int \frac{\partial \phi \cos i \phi}{(1+aa-2a \cos \phi)^{n+1}}$"} \cite{E674}\footnote{The papers \cite{E672}, \cite{E673}, \cite{E674} were essentially part of a single work within the book {\it ``Institutionum calculi integralis"}, Volume 4.} (E674: ``Proof of the extraordinary theorem that had been found by conjecture on the integral formula $\int \frac{\partial \phi \cos i \phi}{(1+aa-2a \cos \phi)^{n+1}}$") and {\it ``Specimen transformationis singularis serierum"} \cite{E710} (E710: ``A specimen of a singular transformation of series"). Whereas the representation in the first paper seems to be less familiar\footnote{It is not listed in the modern table of series and integrals \cite{Zw14} for example.}, the one in the remaining papers can be reduced to the integral:

\begin{equation}
    \label{eq: Laplace Legendre}
    P_n(t) = \dfrac{1}{\pi}\int\limits_{0}^{\pi}\left(t+\sqrt{t^2-1}\cos \varphi\right)^n d\varphi,
\end{equation}
which was given by Laplace  (1749-1827) in his work {\it ``Traité de mécanique céleste"} \cite{La25} (``Treatise on celestial mechanics") and is also attributed to him in the preface to  Volume 16,2 of Series 1 of Euler's {\it Opera Omnia} (\cite{OO162}) . Hence we refer to (\ref{eq: Laplace Legendre}) as the Laplace representation. 

We will consider this representation first (see Section \ref{sec: Euler and the Laplace Representation for the Legendre Polynomials}), before moving on to the less familiar one (see Section \ref{sec: Euler's  other Integral Representation}).

\section{Euler and the Laplace Representation for the Legendre Polynomials}
\label{sec: Euler and the Laplace Representation for the Legendre Polynomials}

In the papers \cite{E672}, \cite{E673}, \cite{E674} and \cite{E710}  Euler considered the family of integrals

\begin{equation}
    \label{eq: Euler Family}
    A_n(a,i):= \int\limits_{0}^{\pi} \dfrac{d \varphi\cos (i \varphi)}{(1+a^2-2a \cos (\varphi))^{n}},
\end{equation}
about which he proved three main results\footnote{In this formula, $i$ denotes an integer and is not to be confused with $\sqrt{-1}$.}. First, he found an explicit formula for the integrals in (\ref{eq: Euler Family}). Secondly, he discovered a functional equation satisfied by the family of integrals. And finally, he stated a difference equation (in $n$) for them.  We will discuss each result separately, stating it and showing how it can be applied to the Legendre polynomials.

The difference equation of interest for us was given by Euler in § 74 of \cite{E673} as follows:

\begin{equation*}
     n(n-1)(1-aa)^2 \int \dfrac{ d \varphi \cos (i \varphi)}{\Delta^{n+1}}
\end{equation*}
\begin{equation*}
%\label{eq: Euler's Difference Equation}
    =(n-1)(2n-1)(1+aa)\int \dfrac{d \varphi \cos(i \varphi)}{\Delta^{n}}+(ii-(n-1)^2)\int \dfrac{d \varphi \cos(i \varphi)}{\Delta^{n-1}},
\end{equation*}
where Euler understood the integrals to be taken from $0$ to $\pi$ and used the abbreviation $\Delta =1-2a\cos (\varphi)+a^2$. 

 Using the notation we introduced in (\ref{eq: Euler Family}), Euler's difference equation reads as follows:

\begin{equation*}
    n(n-1)(1-a^2)^2\cdot A_{n+1}(a,i)
\end{equation*}
\begin{equation*}
    = (n-1)(2n-1)(1+a^2)\cdot A_n(a,i)+(i^2-(n-1)^2)\cdot A_{n-1}(a,i).
\end{equation*}
For $i=0$, defining $A_n(a,0):=A_{n}(a)$ and simplifying the equation, we arrive at:

\begin{equation*}
    n(1-a^2)^2\cdot A_{n+1}(a)=(2n-1)(1+a^2)\cdot A_{n}(a)-(n-1)\cdot A_{n-1}(a).
\end{equation*}
After an index shift $n \mapsto n+1$, this equation reads:

\begin{equation*}
%\label{eq: Intermediate}
    (n+1)(1-a^2)^2\cdot A_{n+2}(a)=(2n+1)(1+a^2)\cdot A_{n+1}(a)-n\cdot A_{n}(a).
\end{equation*}
To simplify this equation even further, we need to reverse Euler's introduction  of the letter $a$ that he made in \cite{E673}. There, he  started from the integral:

\begin{equation*}
    \int \dfrac{d \varphi \cos (i \varphi)}{(\alpha + \beta \cos (\varphi))^n}.
\end{equation*}
The integrals we denoted by $A_n(a,i)$ in (\ref{eq: Euler Family}) arise from this by writing $\alpha = 1+a^2$ and $\beta =-2a$, as Euler also did in § 42 of \cite{E673}. Therefore, the previous difference equation can also be expressed this way: 

\begin{equation*}
    (n+1)(\alpha^2-\beta^2)\cdot A_{n+2}(a) =(2n+1)\alpha \cdot A_{n+1}(a)-n \cdot A_{n}(a).
\end{equation*}
Thus, setting $\alpha =x$ and $\beta= \sqrt{x^2-1}$ and, for the sake of brevity, writing simply $A_n(x)$ instead of $A_n(a)$, we have:

\begin{equation*}
      (n+1)\cdot A_{n+2}(x) =(2n+1)x \cdot A_{n+1}(x)-n \cdot A_{n}(x),
\end{equation*}
which is the difference equation for the Legendre polynomials (\ref{eq: Difference Equation for Legendre}). Moreover, we have $A_{n+1}(x)= \pi P_n(x)$, since by direct calculation:

\begin{equation*}
    A_1(x) = \int\limits_{0}^{\pi} \dfrac{d \varphi}{x+\sqrt{x^2-1}\cos \varphi}= \pi = \pi \cdot P_0(x),
\end{equation*}
provided $x \in [-1,1]$. Further, with the same restrictions on $x$,

\begin{equation*}
    A_2(x) = \int\limits_{0}^{\pi} \dfrac{d \varphi}{(x+\sqrt{x^2-1}\cos \varphi)^2}= \pi \cdot x = \pi \cdot P_1(x). 
\end{equation*}
Therefore, we conclude that
    
    \begin{equation*}
    %\label{eq: Laplace Fake}
        P_n(x) = \dfrac{1}{\pi}\int\limits_{0}^{\pi} \dfrac{d \varphi}{(x+\sqrt{x^2-1}\cos(\varphi))^n},
    \end{equation*}
    which already looks similar to Laplace's representation (\ref{eq: Laplace Legendre}). 
    
To derive the Laplace integral representation for the Legendre polynomials from the last equation, we need a functional equation found by Euler for the integrals $A_n(a,i)$. In § 21 of \cite{E710} Euler proved the formula\footnote{Euler stated the same formula in his papers (\cite{E672}, \cite{E673}, \cite{E674}) without a proof in the first two papers and also remarked himself that the formula was conjectural at that point. The first proof was given by Euler in \cite{E674}.}:

\begin{equation*}
     \binom{n+i}{i} (1-aa)^{-n} \int \Delta^n d \varphi \cos (i\varphi)
\end{equation*}
\begin{equation}
\label{eq: Euler Functional Equation}
    = \binom{-n-1+i}{i}(1-aa)^{n+1} \int \Delta^{-n-1} d \varphi \cos (i \varphi),
\end{equation}
where $\Delta = 1-2a \cos (\varphi)+a^2$ and the integrals are understood to be taken from $0$ to $\pi$. 

 We are again interested in the case $i=0$, in which the binomial coefficients in (\ref{eq: Euler Functional Equation}) become $=1$. After the same substitutions we made earlier to introduce the letter $x$ instead of $a$, in our notation Euler's equation reads:

\begin{equation*}
    A_{-n}(a)= A_{n+1}(a),
\end{equation*}
where we used that by the previous substitutions $1-a^2=\alpha^2-\beta^2=x^2-(x^2-1)=1.$  Combining this with the fact that $A_{n+1}(a)=\pi \cdot P_n(x)$, we also have $A_{-n}=\pi \cdot P_n(x)$. Hence by the definition of $A_{n}(a)$, we also have:

\begin{equation*}
    P_n(x) =\dfrac{1}{\pi}\int\limits_{0}^{\pi} (x+\sqrt{x^2-1}\cos (\varphi))^n d \varphi,
\end{equation*}
which is (\ref{eq: Laplace Legendre}).\\[2mm]
Therefore, using Euler's formulas, we were also able to prove the formula

\begin{equation}
\label{eq: Integral Relation}
   \int\limits_{0}^{\pi} (x+\sqrt{x^2-1}\cos (\varphi))^n d \varphi = \int\limits_{0}^{\pi} (x+\sqrt{x^2-1}\cos (\varphi))^{-n-1} d \varphi, 
\end{equation}
a formula attributed to Jacobi (1854 - 1851) in \cite{OO162}, who actually also proved (\ref{eq: Euler Functional Equation}) in his paper {\it ``Über die Entwickelung des Ausdrucks $(aa-2aa'[\cos \omega \cos \varphi+\sin \omega \sin \varphi \cos(\vartheta-\vartheta')]+a'a')^{-\frac{1}{2}}$"} \cite{Ja43} (``On the expansion of the expression $(aa-2aa'[\cos \omega \cos \varphi+\sin \omega \sin \varphi \cos(\vartheta-\vartheta')]+a'a')^{-\frac{1}{2}}$"). But the argument can be made that Euler could also have proven it, if he had the intention. For, his functional equation (\ref{eq: Euler Functional Equation}) is much more general than (\ref{eq: Integral Relation}). Furthermore, (\ref{eq: Euler Functional Equation}) is a special case of the Euler transformation for the hypergeometric series, i.e., the formula

\begin{equation*}
    _2F_1(a,b,c;x)= (1-x)^{c-a-b}{}_2F_1(c-a, c-b,c; x),
\end{equation*}
which Euler proved in § 9 of \cite{E710}. Here, $_2F_1$ is the Gaussian hypergeometric series. Therefore,  (\ref{eq: Integral Relation}) can be understood as a special case of a certain property of the Gaussian hypergeometric function.

For the sake of completeness, we also state an explicit formula that Euler gave for the integrals $A_n(a,i)$. The formula reads as follows (see, e.g., § 43 in \cite{E673}):

\begin{equation*}
    \int\limits_{0}^{\pi} \dfrac{d \varphi \cos (i \varphi)}{(1+aa-2a \cos (\varphi))^{n+1}}= \dfrac{\pi a^i}{(1-a^2)^{2n+1}} \cdot V,
\end{equation*}
where

\begin{equation*}
    V= \binom{n+i}{i}+\binom{n+i}{i+1}\binom{n-i}{1}a^2+\binom{n+i}{i+2}\binom{n-i}{2}a^4+\cdots
\end{equation*}
While Euler himself noted in \cite{E672} and \cite{E673} that the formula for $V$ was conjectural, he gave a proof in \cite{E674}.

\section{Euler's other Integral Representation}
\label{sec: Euler's  other Integral Representation}

As mentioned in the introduction, the other integral representation for the Legendre polynomials originates in Euler's paper \cite{E606}, a paper originally devoted to continued fractions. Having shown in his papers {\it ``De fractionibus continuis observationes"} \cite{E123} (E123: ``Observations on continued fractions")  and {\it ``Methodus inveniendi formulas integrales, quae certis casibus datam inter se teneant rationem, ubi simul methodus traditur fractiones continuas summandi"} \cite{E594} (E594: ``A method to find integral formulas which in certain case are in a given ratio to each other, where at the same a method is given to sum continued fractions")  how homogeneous difference equations with linear coefficients can be solved by means of integrals, in  \cite{E606}, Euler applied his theory to the family of integrals 

\begin{equation*}
%\label{eq: Euler Integral Difference Equation}
    G(n):=\int\limits_{\frac{b-\sqrt{b^2-a^2c}}{c}}^{\frac{b+\sqrt{b^2-a^2c}}{c}} \dfrac{x^ndx}{\sqrt{a^2-2bx+cx^2}}
\end{equation*}
and showed that those integrals satisfy the equation:

\begin{equation*}
    na^2\cdot G(n-1)= (2n+1)b\cdot G(n)-(n+1)c\cdot G(n+1).
\end{equation*}
Specializing to the  case $a=c=1$ and writing $t$ for $b$, this equation reduces to the difference equation for the Legendre polynomials (\ref{eq: Difference Equation for Legendre}) such that we can conclude

\begin{equation*}
C(t)\cdot    P_n(t) =   \int\limits_{t-\sqrt{t^2-1}}^{t+\sqrt{t^2-1}}\dfrac{x^n dx}{\sqrt{1-2xt+x^2}}.
\end{equation*}
$C(t)$ is a function independent of $n$ that we have to determine. 

To do this, we evaluate the integrals for $n=0$ and $n=1$ explicitly. By a formal calculation, we have:

\begin{equation*}
    \int\limits_{t-\sqrt{t^2-1}}^{t+\sqrt{t^2-1}}\dfrac{x^0 dx}{\sqrt{1-2xt+x^2}} = \log(-1)
\end{equation*}
and

\begin{equation*}
    \int\limits_{t-\sqrt{t^2-1}}^{t+\sqrt{t^2-1}}\dfrac{x^1 dx}{\sqrt{1-2xt+x^2}}= \log(-1)\cdot t.
\end{equation*}
Hence, it turns out that the function $C(t)$ is actually a constant. Unfortunately, $\log(-1)$ is not uniquely defined; we have to choose a specific branch of the complex logarithm. Thus, if we  write

\begin{equation}
   \label{eq: New Legendre Integral}
    P_n(t)= \dfrac{1}{\log(-1)} \cdot \int\limits_{t-\sqrt{t^2-1}}^{t+\sqrt{t^2-1}}\dfrac{x^n dx}{\sqrt{1-2xt+x^2}},
\end{equation}
this equation gives us the correct formula, as long as we do not care about the ambiguous meaning of the symbol $\log(-1)$. If we want to avoid this ambiguity, we can use the indefinite integral

\begin{equation*}
%\label{eq: Claim}
    \int \dfrac{x^ndx}{\sqrt{1-2xt+x^2}}
\end{equation*}
\begin{equation*}
    = P_n(t)\operatorname{artanh}\left(\dfrac{x-t}{\sqrt{1-2xt+x^2}}\right)+ Q_n(x,t) \sqrt{1-2xt+x^2}+C
\end{equation*}
where $Q_n$ is a polynomial in $x$ and $t$ and $C$ is a constant of integration. This formula can be proved by induction and the difference equation for the Legendre polynomials (\ref{eq: Difference Equation for Legendre}). Note that we do not need to know the polynomial $Q_n$ explicitly, since from the boundaries of integration in the previous equation we have $Q_n(x,t)\sqrt{1-2xt+x^2}=0$.

\section{Another Approach leading to similar Insights}
\label{sec: Connection of this Result to similar Insights}

As seen in the introduction, the Legendre polynomials can also defined as coefficients of a generating function (see (\ref{eq: Definition Legendre polynomials})). In this section we intend to show that Euler was led to the aforementioned generating function and  (\ref{eq: Difference Equation for Legendre}) in his investigations in the papers {\it ``Observationes analyticae"} \cite{E326} (E326: ``Analytical Observations") and {\em ``Varia artificia in serierum indolem inquirendi"} \cite{E551} (E551: ``Various Artifices to investigate the Nature of Series"). 
In both these papers, Euler considered the general trinomials $(a+bx+cx^2)^n$ for $n \in \mathbb{N}$ and $a,b,c \in \mathbb{C}$ and  was interested in the sequence of the coefficients of the power $x^n$ in the expansion of those trinomials. Therefore, let us assume

\begin{equation*}
%\label{eq: Assumption}
    (a+bx+cx^2)^n = a^n +\cdots + B_n\cdot x^n+ \cdots +c^nx^{2n},
\end{equation*}
where the coefficient $B_n$ depends on the numbers $a$, $b$ and $c$. From these coefficients Euler then defined a Taylor series:

\begin{equation}
\label{eq: Taylor Series}
    F(x):= \sum_{n=0}^{\infty} B_n x^n.
\end{equation}
The problem he solved in \cite{E326} and \cite{E551} is\footnote{The case $a=b=c=1$ gives coefficients comprising sequence A002426 in the Online Encyclopedia of Integer Sequences (\url{https://oeis.org/}).}: What is a closed expression for $F(x)$? His answer to the problem (see, e.g., § 21 of \cite{E551}) reads:
\begin{equation*}
%\label{eq: General Solution}
    F(x):= \dfrac{1}{\sqrt{1-2bx+(b^2-4ac)x^2}}.
\end{equation*}
This expression bears close resemblance to the generating function for the Legendre polynomials in (\ref{eq: Definition Legendre polynomials}). Indeed, writing $t$ for $b$ and setting $b^2-4ac=1$ in Euler's solution, we arrive at the generating functions for the Legendre polynomials. 

Another of Euler's important contributions in \cite{E551} was his discovery of  the difference equation between three consecutive Legendre polynomials, i.e, (\ref{eq: Difference Equation for Legendre}). Euler gave  the following relation among the consecutive expansion coefficients of the Taylor series for $F$ in (\ref{eq: Taylor Series}) in § 22 of  \cite{E551}:

\begin{equation*}
%\label{eq: Euler Difference Legendre}
    r= bq+\dfrac{n-1}{n}(bq-p).
\end{equation*}
Here, Euler suppressed the dependence of $r,p,q$ on $n$ and also on the letter $b$, which is the same $b$ as in explicit expression for the function $F$. But Euler told us that, $r$ being the $n$-th term in the sequence of coefficients, we obtain $q$ from $r$ by writing $n-1$ instead of $n$ in $r$ and likewise, $p$ from $r$ by writing $n-2$ instead of $n$ in $r$. Thus, writing $t$ instead of $b$ and presenting everything in modern fashion with indices, i.e., setting $r=B_{n}(t)$, $q=B_{n-1}(t)$ and $p=B_{n-2}(t)$, Euler's difference equation becomes:

\begin{equation*}
    B_n(t)= t \cdot B_{n-1}(t)+\dfrac{n-1}{n}(t \cdot B_{n-1}(t)-B_{n-2}(t)).
\end{equation*}
A simple calculation shows that this equation is equivalent to:

\begin{equation*}
    n B_n(t)=(2n-1)t\cdot B_{n-1}(t)-(n-1)B_{n-2}(t).
\end{equation*}
After a shift $n \mapsto n+1$, this equation is to be seen to be identical to  the difference equation for the Legendre polynomials (\ref{eq: Difference Equation for Legendre}). 

As a side note we mention that, we can use Euler's approach and define the Legendre polynomials as follows:

\begin{Def}
The Legendre polynomial $P_n(t)$ is given as the coefficient of the power $x^n$ the expansion of the the following trinomial:

\begin{equation*}
   T(t,x)= \left(\dfrac{t-1}{2}+xt+\dfrac{t+1}{2}x^2\right)^n.
\end{equation*}
\end{Def}
This definition allows to calculate the Legendre polynomials, since, as we saw above,  for this particular choice of the letters $a,b,c$ in the assumed form of the expanded trinomial the function $F$ in (\ref{eq: Taylor Series}) reduces to the generating function in (\ref{eq: Definition Legendre polynomials}). 

\section{Conclusion}
\label{sec: Conclusion}

In this note we showed that Euler's discoveries lead to two integral representations for the Legendre polynomials. The first equation (\ref{eq: Laplace Legendre}) that we derived from Euler’s many investigations (\cite{E672}, \cite{E673}, \cite{E674} and \cite{E710}) is familiar and is attributed to Laplace (\cite{OO162}), whereas the other equation (\ref{eq: New Legendre Integral}),  which follows from Euler's formulas in (\cite{E606}), seems to be less known, since it does not appear in a standard integral reference tables (\cite{Zw14}). Furthermore, we showed that Euler's investigations in \cite{E326} and \cite{E551} contain the generating function for the Legendre polynomials (\ref{eq: Definition Legendre polynomials}) as a special case and the difference equation (\ref{eq: Difference Equation for Legendre}) was stated explicitly in \cite{E551}. Additionally, the papers \cite{E606}, \cite{E673} include it as a special case, the general equation being (\ref{eq: Euler Functional Equation}). Unfortunately, it seems that Euler himself never made the connection between his findings, although they all led him to (\ref{eq: Difference Equation for Legendre}).

Euler himself did not know the Legendre polynomials and their definition; they were introduced by Legendre in \cite{Le85} in his studies of the gravitational potential after Euler's death. Therefore, it is even more interesting that they also appeared in Euler's investigations, which were not motivated by a problem in physics. Indeed, Euler's papers rather seem like curious intellectual challenges without any clear origin. Hence this should provide us with a lot of motivation to consider more of Euler's papers from a modern perspective. They might even lead to new interpretations of familiar results and contain more insights attributed to other mathematicians like (\ref{eq: Euler Functional Equation}), a formula attributed to Jacobi in \cite{OO162} who actually proved it later than Euler in \cite{Ja43}.

\end{document}